\documentclass[a4paper,12pt]{article}
\usepackage{latexsym}
\usepackage{amssymb}
\usepackage[mathscr]{eucal}
\newcommand{\rr}{\mathbb{R}}

\newcommand{\nn}{\mathbb{N}}
\newcommand{\zz}{\mathbb{Z}}
\newtheorem{Theorem}{Theorem}

\newtheorem{Lemma}[Theorem]{Lemma}
\newtheorem{Corollary}[Theorem]{Corollary}

\begin{document}
\title{Convergence of ray sequences of Pad\'e approximants for
$\phantom{}_2F_1(a,1;c;z)$, for $c>a>0$.}
\author{K. Driver\thanks{John Knopfmacher Centre for Applicable Analysis and
Number Theory, School of Mathematics, University of the
Witwatersrand, Johannesburg, South Africa} \ and K. Jordaan
\thanks{Department of Mathematics and Applied Mathematics,
University of Pretoria, Pretoria, 0002, South Africa}}
\date{}
\maketitle

\bigskip

\begin{abstract}
The Pad\'e table of $\phantom{}_2F_1(a,1;c;z)$ is normal for
$c>a>0$ (cf. \cite{3}). For $m \geq n-1$ and $c \notin
{\zz}^{\phantom{}^-}$, the denominator polynomial $Q_{mn}(z)$ in
the $[m/n]$ Pad\'e approximant $P_{mn}(z)/Q_{mn}(z)$ for
$\phantom{}_2F_1(a,1;c;z)$ and the remainder term
$Q_{mn}(z)\phantom{}_2F_1(a,1;c;z)-P_{mn}(z)$ were explicitly
evaluated by Pad\'e (cf. \cite{2}, \cite{5} or \cite{7}). We show
that for $c>a>0$ and $m\geq n-1$, the poles of
$P_{mn}(z)/Q_{mn}(z)$ lie on the cut $(1,\infty)$. We deduce that
the sequence of approximants $P_{mn}(z)/Q_{mn}(z)$ converges to
$\phantom{}_2F_1(a,1;c;z)$ as $m \to \infty$, $ n/m \to \rho$ with
$0<\rho \leq 1$, uniformly on compact subsets of the unit disc
$|z|<1$ for $c>a>0$.
\end{abstract}
\bigskip

\noindent AMS MOS Classification:\quad 41A21, 30E15 \bigskip


\newpage
\section{The Pad$\mathbf{\acute{e}}$ approximant for \boldmath{$_2F_1(a,1;c;z)$} \\when
{\boldmath{$m \geq n-1$}}.}

The Gauss hypergeometric function, or $\phantom{}_2F_1$, is
defined by
\[\phantom{}_2F_1(a,b;c;z)=
1+\sum_{k=1}^{\infty}\frac{(a)_k(b)_kz^k}{(c)_kk!},\quad |z|<1,\]
where the parameters $a, b, c$ and $z$ may be real or complex and
\[(\alpha)_k=\alpha(\alpha+1)\ldots(\alpha+k-1)
=\frac{\Gamma(\alpha+k)}{\Gamma(\alpha)}\] is Pochhammer's symbol.
When $a=-n$ is a negative integer, the series terminates and
reduces to a polynomial of degree $n$. \bigskip

The $[m/n]$ Pad\'e approximant for a formal power series
\[f(z)=\sum_{k=0}^{\infty}t_kz^k\]
is a rational function $P_{mn}(z)/Q_{mn}(z)$ of type $(m,n)$ such
that
\begin{equation}\label{1}
f(z)Q_{mn}(z)-P_{mn}(z)=0(z^{m+n+1}) \mbox{ as }z \rightarrow 0.
\end{equation} The $[m/n]$ Pad\'e approximants for $f$ can be
arranged to form the Pad\'e table of $f$ (cf. \cite{2}). A Pad\'e
approximant is normal if it occurs only once in the Pad\'e table
and the Pad\'e table is normal if each entry in the table is
normal. The following results have been proved for the Pad\'e
approximants for $\phantom{}_2F_1(a,1;c;z)$.
\bigskip

\begin{Theorem}[cf. \cite{3}, Corollary 6.1]\label{normality}
The Pad\'e table for the hypergeometric series
$\phantom{}_2F_1(a,1;c;z)$ with $c>a>0$ is
normal.\end{Theorem}\bigskip

\begin{Theorem}[cf. \cite{5}, p. 389, \cite{2},  \cite{7}]
Let $c \notin {\zz}^{\phantom{}^-}$ and let $m \geq n-1$. Then the
denominator polynomial in the $[m/n]$ Pad\'e approximant
$P_{mn}(z)/Q_{mn}(z)$ for $_2F_1(a,1;c;z)$ is given by
\begin{equation}\label{Q}
Q_{mn}(z)=\phantom{}_2F_1(-n,-a-m;-c-m-n+1;z)\end{equation} and
\begin{eqnarray}\label{remainder}
R_{mn}(z) & = & Q_{mn}(z)\phantom{1}_2F_1(a,1;c;z)-P_{mn}(z)\nonumber \\
& = &
S_{mn}\phantom{1}z^{m+n+1}\phantom{1}_2F_1(a+m+1,n+1;c+m+n+1;z)
\end{eqnarray}
where \begin{equation}\label{S}
S_{mn}=n!\frac{(a)_{m+1}(c-a)_n}{(c)_{m+n}(c+m)_{n+1}}\phantom{1}.
\end{equation}
\end{Theorem}\bigskip

\begin{Theorem}[cf. \cite{7}, Theorem 6.2]
If $a, c, c-a \notin {\zz}^{\phantom{}^-}$, the Pad\'e
approximants for $_2F_1(a,1;c;z)$ are normal for $m \geq n-1$.
\end{Theorem}\bigskip

The numerator polynomial $P_{mn}(z)$ is determined by (\ref{1})
since $Q_{mn}(z)$ is known. Clearly, $P_{mn}(z)$ is the polynomial
we obtain from the first $(m+1)$ terms in the product
\[
_2F_1(a,1;c;z) \phantom{1}_2F_1(-n,-a-m;-c-m-n+1;z).\] Thus
\[P_{mn}(z)=\sum_{r=0}^{m}p_rz^r\] where for $ 0 \leq r \leq m$,
\[p_r = \sum_{l=0}^{r}\frac{( a)_{r-l}(-n)_l(-a-m)_l}{(-c-m-n+1)_l
(c)_{r-l}l!}.\]
$P_{mn}(z)$ is not, in general, a $\phantom{}_2F_1$ hypergeometric polynomial.
For example, an elementary calculation shows that for $a=2$, $c=6$, $m=3$ and $n=4$,
\[P_{34}(z) = 1-\frac{4}{3}z+\frac{344}{693}z^2-\frac{1}{22}z^3\]
which is not equal to $\phantom{}_2F_1(-3,\alpha;\beta;z)$ for any
$\alpha, \beta$. Taking a diagonal entry from the Pad\'e table,
for instance, letting $a=3.2$, $c=5.44$, $m=3$ and $n=3$ results
in
\[P_{33}(z) = 1-1.19337z+0.317021z^2-0.000851604z^3\]
which also is not equal to $\phantom{}_2F_1(-3,\alpha;\beta;z)$
for any $\alpha, \beta$. \bigskip

The location and behavior of the zeros and poles of Pad\'e
approximants for various special functions, as well as the
asymptotic zero and pole distribution, has been studied by many
authors, most notably E. Saff and R. Varga (cf. \cite{6}). We
shall use the Rodrigues' formula for $\phantom{}_2F_1(-n,b;d;z)$
to establish the location of the zeros of $Q_{mn}(z)$ in various
cases. In the case when the Pad\'e table of
$\phantom{}_2F_1(a,1;c;z)$ is normal, ie $c>a>0$, we will show
that the poles of the approximant lie in the interval
$(1,\infty)$.
\bigskip

The convergence of different types of sequences in the Pad\'e
table has been studied extensively. In \cite{Perron}, O. Perron
studies convergence in the Pad\'e table for the exponential
function and M.G. de Bruin \cite{debruin} shows that the
convergence behavior in the Pad\'e table for $_1F_1(1;c;z)$ with
$c\notin {\zz}^{\phantom{}^-}$ is similar to that of
$e^z$.\bigskip

Using the explicit form of the remainder $R_{mn}(z)$ given in
(\ref{remainder}), we show that the sequence of Pad\'e
approximants for $_2F_1(a,1;c;z)$, $m\geq n-1$ converges to
$_2F_1(a,1;c;z)$ as $m\to \infty$, $ n/m \to \rho$ with $0<\rho
\leq 1$ on compact subsets of $|z|<1$ for $c>a>0$.
\bigskip

\section{The zeros of \boldmath{$_2F_1(-n,b;d;z)$} and \boldmath{$Q_{mn}(z)$}.}
Rodrigues' formula states (cf. \cite{1}, p. 99)
\[z^{d-1}(1-z)^{b-d-n}\phantom{}_2F_1(-n,b;d;z)=\frac{1}{(d)_n}\frac{d^n}{dz^n}
[z^{d-1+n}(1-z)^{b-d}].\] Let $g_{\phantom{}_{l}}(z)$ be an
arbitrary polynomial of degree $l<n$. Then
\begin{eqnarray}\label{3} \lefteqn{(d)_n\int_p^q
z^{d-1}(1-z)^{b-d-n}\phantom{}_2F_1(-n,b;d;z)g_{\phantom{}_l}(z)dz}\nonumber\\
&=&\int_p^q\{D^n[z^{d-1+n}(1-z)^{b-d}]\}g_{\phantom{}_l}(z)dz.\end{eqnarray}
Integrating the right hand side of (\ref{3}) by parts $n$ times,
each time differentiating $g_{\phantom{}_l}(z)$ and integrating
the expression in curly brackets, we obtain
\begin{eqnarray}\label{4}
\lefteqn{\left[\sum_{k=1}^{n}{(-1)^{k-1}D^{n-k}[z^{d-1+n}(1-z)^{b-d}]}
D^{k-1}[g_{\phantom{}_l}(z)]\right]_p^q}\nonumber\\
&+&(-1)^n\int_p^q{z^{d-1+n}(1-z)^{b-d}}D^n[g_{\phantom{}_l}(z)]dz.\end{eqnarray}
Each term in the sum in (\ref{4}) above contains a product of
powers of $z$ and powers of $(1-z)$, where the lowest and highest
powers of $z$ are $d$ and $(d+n-1)$ respectively. The lowest and
highest powers of $(1-z)$ are $(b-d-n+1)$ and $(b-d)$
respectively. Moreover, each term occurring in the sum in
(\ref{4}) is $O(z^{b+l})$ as $z \rightarrow \pm \infty$, where $l
\leq n-1$.
\begin{Theorem}\label{zeros} Let $F(z)=F(-n,b;d;z)$, where $n \in \nn$, $b, d \in \rr$
and $d \notin {\zz}^{\phantom{}^-}$.
\begin{itemize}
\item[(i)]For $d>0$ and $b>d+n-1$, all $n$ zeros of
$F$ are real and simple and lie in $(0,1)$.
\item[(ii)]For $b<1-n$ and $d<b+1-n$, all $n$ zeros of $F$ are
real and simple and lie in $(1,\infty)$.
\item[(iii)]For $b<1-n$ and $d>0$, all zeros of $F$ are real and
simple and lie in $(-\infty,0)$.
\end{itemize} \end{Theorem}\bigskip

\noindent{\it Proof. }
\begin{itemize}
\item[(i)]Putting $p=0$ and $q=1$ in (\ref{3}) and using
(\ref{4}), we have \begin{equation}\label{5} \int_0^1
z^{d-1}(1-z)^{b-d-n}F(z)g_{\phantom{}_l}(z)dz=0\end{equation}
since $d>0$, $b>d+n-1$ and $D^n[g_{\phantom{}_l}(z)]=0$. The
function \\
$z^{d-1}(1-z)^{b-d-n}$ does not change sign in $(0,1)$, so from
(\ref{5}) with $l=0$, we know that $F(z)$ changes sign at least
once in $(0,1)$. Let $z_1, \ldots ,z_l$ be the zeros of $F$ in
$(0,1)$ of odd multiplicity. If $l<n$ we can let
\[g_{\phantom{}_l}(z)=\prod_{r=1}^l (z-z_r)\] and obtain a contradiction from
(\ref{5}) since $F(z)g_{\phantom{}_l}(z)$ does not change sign in
$(0,1)$. Therefore $l=n$ and $F(z)$ has $n$ zeros of odd
multiplicity in $(0,1)$ and the result follows. \item[(ii)]
Putting $p=1$ and $q=\infty$ in (\ref{3}) and using (\ref{4}), we
see that for $b<1-n$ and $d<b+1-n$, \[ \int_1^{\infty}
z^{d-1}(1-z)^{b-d-n}F(z)g_{\phantom{}_l}(z)dz=0\] A similar
argument to that used in (i) completes the proof. Note that
$z^{d-1}(1-z)^{b-d-n}$ does not change sign for $z>1$.
\item[(iii)] Putting $p=-\infty$ and $q=0$ in (\ref{3}) and using
(\ref{4}), we have that if $b<1-n$ and $d>0$, \[\int_{-\infty}^0
z^{d-1}(1-z)^{b-d-n}F(z)g_{\phantom{}_l}(z)dz=0\] and again,
$z^{d-1}(1-z)^{b-d-n}$ does not change sign for $z<0$. \quad
\end{itemize}
\bigskip

{\bf Remark:} The results proved in Theorem \ref{zeros} were known
to Klein \cite{4} but his proofs are somewhat less transparent
than the simple use of Rodrigues' formula.\bigskip

\begin{Corollary}\label{poles} For $c \notin \zz^{\phantom{}^-}$ and $m \geq n-1$,
the poles of the $[m/n]$ Pad\'e approximant for
$\phantom{}_2F_1(a,1;c;z)$ lie in the intervals
\begin{itemize}
\item[(i)] $(0,1)$ if $a<c<1-m-n$ \item[(ii)]$(1,\infty)$ if
$c>a>n-m-1$ \item[(iii)]$(-\infty,0)$ if $a>n-m-1$ and $c<1-m-n$.
\end{itemize}
\end{Corollary}\bigskip

\noindent{\it Proof. }Let $b=-a-m$ and $d=-c-m-n+1$ in Theorem
\ref{zeros}. \quad $\Box$ \bigskip

\section{Convergence in the Pad$\mathbf{\acute{e}}$ table for \\
\boldmath{$_2F_1(a,1;c;z)$}.}

In this section we prove convergence results for the Pad\'e
approximants of $_2F_1(a,1;c;z)$. \bigskip

\begin{Lemma}\label{remainderlemma}
For $m\geq{n-1}$ and $c>a>0$ we have that
\[R_{mn}(z)=Q_{mn}(z)\phantom{}_2F_1(a,1;c;z)-P_{mn}(z)\] tends
to zero uniformly in $z$ as $m \to \infty$, $n/m \to \rho$ with
$0<\rho \leq 1$ on compact subsets of $|z|<1$.\end{Lemma}\bigskip

\noindent{\it Proof. } For $c>a>0$ we have
\begin{eqnarray}\label{upperbound}\lefteqn{\left|
\phantom{}_2F_1\left(a+m+1,n+1;c+m+n+1;z\right) \right|}\nonumber
\\
& \leq & \sum_{k=0}^{\infty}\frac{(a+m+1)_k(n+1)_k}{(c+m+n+1)_k}
\frac{|z|^k}{k!}.\end{eqnarray} For $c-a-1>0$, the series on the
right converges at $z=1$ and we have for $|z|<1$ (cf. \cite{1},
Theorem 2.2.2, p. 66),
\begin{eqnarray}\label{bound1}
\lefteqn{\left|\phantom{}_2F_1\left(a+m+1,n+1;c+m+n+1;z\right)
\right|}\nonumber \\& \leq &
\frac{\Gamma(c+m+n+1)\Gamma(c-a-1)}{\Gamma(c+m)
\Gamma(c-a+n)}.\end{eqnarray} Then for $m\geq n-1$ , $c-1>a>0$ and
$|z|<1$ we have the estimate from (\ref{remainder}) , (\ref{S})
and (\ref{bound1})
\begin{eqnarray*}
\left|R_{mn}(z)\right|& \leq &
k\frac{n!(a+1)_m}{\Gamma(c+m+n}|z|^{m+n+1}\\
& = &k'\frac{n!(a+1)_m}{(c)_n(c+n)_m}|z|^{m+n+1}\\
& \rightarrow & 0 \mbox{ as } m\rightarrow \infty,\end{eqnarray*}
since $c>1$ and $c+n>a+1>1$ implies that $\frac{n!}{(c)_n}<1$ and
$\frac{(a+1)_m}{(c+n)_m}<1$ for all $m,n \in \nn$. The constants
$k$ and $k'$ are independent of $n$ and $m$. For $c-a-1<0$ the
series on the right of (\ref{upperbound}) diverges at $z=1$ but it
follows from \cite{1}, p.63, Theorem 2.1.3, for $|z|<1$, that
\begin{eqnarray}\label{bound2}
\lefteqn{\left|\phantom{}_2F_1\left(a+m+1,n+1;c+m+n+1;z\right)
\right|}\nonumber \\ & \leq & |(1-z)|^{c-a-1}\frac{\Gamma(c+m+n+1)
\Gamma(c-a-1)}{\Gamma(n+1)\Gamma(a+m+1)}.\end{eqnarray} Then for
$m\geq n-1$, $c-a-1<0$ and $|z|<1$,  we have the estimate from
(\ref{remainder}), (\ref{S}) and (\ref{bound2}),
\begin{eqnarray*}
\lefteqn{\left|R_{mn}(z)\right|}\\& \leq &
\frac{n!(a)_{m+1}(c-a)_n\Gamma(c+m+n+1)\Gamma(a-c+1)}
{(c)_{m+n}(c+m)_{n+1}\Gamma(n+1)\Gamma(a+m+1)}|z|^{m+n+1}|1-z|^{c-a-1}\\
& \leq &
w\frac{(c-a)_n}{(c+m)_n}|z|^{m+n+1}|1-z|^{c-a-1}\\
& \rightarrow & 0 \mbox{ as } m \rightarrow \infty,
\end{eqnarray*} since $0<c-a<c+m$ implies $\frac{(c-a)_n}{(c+m)_n}<1$ for
all $m,n \in \nn$. The constant $w$ is independent of $m$ and
$n$.$\Box$\bigskip

\begin{Theorem} Let $a, c, c-a \notin {\zz}^{\phantom{}^-}$ and $m \geq
n-1$. The sequence of Pad\'e approximants $P_{mn}(z)/Q_{mn}(z)$
converges to $_2F_1(a,1;c;z)$ for $m\to \infty$, $ n/m \to \rho$
with $0<\rho \leq 1$, uniformly in $z$ on compact subsets of the
unit disc $|z|<1$ for $c>a>0$.\end{Theorem}\bigskip

\noindent{\it Proof. } We know from Corollary \ref{poles} that all
the zeros of $Q_{mn}(z)$ lie on the cut $(1,\infty)$ for $c>a>0$,
$m\geq n-1$. Therefore we can divide the result of Lemma
\ref{remainderlemma} by $Q_{mn}(z)$ since
$\left|Q_{mn}(z)\right|>0$ for all $z$ in the unit disc.$\Box$
\bigskip

{\bf Acknowledgement:} The authors are indebted to the referee who
suggested significant improvements.

\bigskip

\begin{tabbing}
e-mail addresses: \= kathy@maths.wits.ac.za \\
\> kjordaan@scientia.up.ac.za
\end{tabbing}

\end{document}